\documentclass[twoside,10pt,reqno]{amsart}

\usepackage{upref,amsxtra,amssymb,amscd}
\usepackage{varioref}
\usepackage{verbatim}
\usepackage{epsfig}

\usepackage{epic,eepic,eucal}
\usepackage{enumerate}

\def\E{          \mathcal E}
\def\D{          \mathcal D}

\def\B{          \mathcal B}

\def\A{          \mathcal A}

\def\a{         \alpha}

\newcommand{\NN}{{\mathbb N}}
\newcommand{\RR}{{\mathbb R}}
\newcommand{\TT}{{\mathbb T}}

\newcommand{\ZZ}{{\mathbb Z}}

\newcommand{\QQ}{{\mathbb Q}}

\newcommand{\DR}{{\mathcal DR}}

\newtheorem{theo}{\sc Theorem}[section]
\newtheorem{prop}[theo]{\sc Proposition}
\newtheorem{lemma}[theo]{\sc Lemma}
\newtheorem{sublemm}[theo]{\sc Sublemma}
\newtheorem{lemm}[theo]{\sc Lemma}
\newtheorem{coro}[theo]{\sc Corollary}

\theoremstyle{definition}

\newtheorem{defn}[theo]{\sc Definition}

\theoremstyle{remark}

\newtheorem{remark}[theo]{\sc Remark}

\numberwithin{equation}{section}

\title[Rotations and the Shrinking Target Properties]
      {On Circle Rotations and the Shrinking Target Properties }

\author[Jimmy Tseng]{}

\subjclass{Primary: 37A05, 37E10; Secondary: 11K60}
 \keywords{Recurrence, circle rotations, logarithm laws, continued fractions, Diophantine approximation}
 
\email{jtseng@brandeis.edu}

\begin{document}
\maketitle

\centerline{\scshape Jimmy Tseng }
\medskip
{\footnotesize
 \centerline{Department of Mathematics, Brandeis University}
   \centerline{Waltham, MA 02454, USA}
} 

\bigskip

 %\centerline{(Communicated by the associate editor name)}
 
\begin{abstract} We generalize the monotone shrinking target property~(MSTP) to the $s$-exponent monotone shrinking target property~($s$MSTP) and give a necessary and sufficient condition for a circle rotation to have $s$MSTP.

Using another variant of MSTP, we obtain a new, very short, proof of a known result, which concerns the behavior of irrational rotations and implies a logarithm law similar to D. Sullivan's logarithm law for geodesics.

\end{abstract}

\section{Introduction}

Let $(M, \mu, T)$ be a measure preserving dynamical system with $\mu(M)$ finite.  The Poincar\'e recurrence theorem asserts that for any fixed measurable set $A$ and for almost every $x \in A$, $T^n(x) \in A$ for infinitely many $n \in \NN$.  If we replace $A$ with a sequence of measurable sets $\{A_n\}$, we can still ask for the measure of the limsup set: $$\{x \in M \mid T^n(x) \in A_n \textrm{ for infinitely many } n \in \NN \}= \textrm{limsup}T^{-n}A_n.$$  

If $\sum \mu(A_n) < \infty$, the convergence case of the Borel-Cantelli lemma implies that the measure of the limsup set is zero.  Otherwise, the situation is more complicated and more interesting.  A general definition, found in~\cite{KM}, applies to this situation:  

\begin{defn}\label{bc} A sequence of measurable sets ${\lbrace  A_n \rbrace}_{n \in \NN }$, such that 
\begin{eqnarray} \label{div} \sum_{n=1}^{\infty} \mu (A_n) =  \infty, \end{eqnarray}
is called a {\it Borel-Cantelli (BC) sequence} for $T$ if 
$$\mu \left(\lim \sup T^{-n}A_n\right)=\mu(M).$$
\end{defn}

\noindent Unfortunately, the divergence case of the Borel-Cantelli lemma is not usually helpful for finding BC sequences since this case of the lemma requires independent sets, a condition that almost never occurs for dynamical systems.  To obtain a BC sequence, we must impose restrictions.  If, for a particular dynamical system, all sequences of measurable sets $\{A_n\}$ that satisfy (\ref{div}) and certain restrictions are BC, we obtain what is called a dynamical Borel-Cantelli lemma.  The earliest example known to the author of such a lemma, Theorem~\ref{theoKurzweil}, in which only sequences of balls centered at a fixed point and with weakly monotonically decreasing radii are allowed, was proved in 1955 by J.~Kurzweil~\cite{Ku}.  

Dynamical Borel-Cantelli lemmas can have applications beyond dynamical systems.  For example, D.~Kleinbock and G.~Margulis discovered such a lemma on homogeneous spaces in which only sequences of complements of balls centered at a fixed point are allowed.  Their lemma and its generalization have important implications for geometry and number theory~\cite{KM}.  One of these implications is a sharper version of D.~Sullivan's logarithm law for geodesics, which is similar to our Corollary~\ref{corLogLaw}.  A more comprehensive introduction to dynamical Borel-Cantelli lemmas, including other examples, can be found in~\cite{CK} and~\cite{KM}.

On a metric space, a sensible restriction to impose would be to consider sequences of balls centered at a fixed point.  Given such sequences, a dynamical Borel-Cantelli lemma allows us to find a set of full measure whose elements return infinitely many times to balls centered at a point whenever the sum of the measures of these balls diverges.  The balls in such a sequence are called shrinking targets by R.~Hill and S.~Velani~\cite{HV}.  Keeping with their terminology, let us say that systems for which we can show a dynamical Borel-Cantelli lemma using sequences of balls centered at a fixed point have a shrinking target property.  Precise definitions and examples follow.  

\subsection{Definitions and Examples.}  Let $(M, \mu, T)$ (or briefly $T$) denote a measure preserving dynamical system where $M$ is a metric space and $\mu$ is a measure on $M$.  Let ${\mathcal B}_M$ denote the Borel $\sigma$-algebra of $M$.  Require that $\mu(M) < \infty$ and that the domain of $\mu \supseteq {\mathcal B}_{M}$.

Define a \textbf{radius sequence} to be a function $r:\NN \rightarrow \RR_{\geq 0}$.  Let $r_n := r(n)$ and denote the set of radius sequences by ${\mathcal R}$.  Define an \textbf{admissible set of radius sequences}, ${\mathcal A}$, to be a non-empty subset of ${\mathcal R}$.

Let us begin with a general definition from which we will specify distinguished special cases:

\begin{defn}[$\A$STP]
The dynamical system $(M, \mu, T)$ has the \textbf{$\A$-shrinking target property ($\A$-STP)} if, for any $x \in M$ and any $r \in {\mathcal A}$ such that $A_n := B(x, r_n)$ satisfies (\ref{div}), $\{A_n \}$ is BC for $T$.
\end{defn}

\noindent In the literature, ${\mathcal R}$-STP is known as the \textbf{shrinking target property (STP)}~\cite{Fa}.  Denote $${\mathcal DR} := \{r \in {\mathcal R} \mid r_n \geq r_{n+1} \textrm{ for all } n \in \NN\}.$$  Likewise, ${\mathcal DR}$-STP is known in the literature as the \textbf{monotone shrinking target property (MSTP)}~\cite{Fa}. Note that STP $\Rightarrow$ MSTP.\bigskip

Expanding maps of the circle~\cite{Ph} and Anosov diffeomorphisms~\cite{Do}, for example, have STP.  Toral translations, for example, do not~\cite{Fa}.  However, whether a toral translation has MSTP or not depends on a Diophantine condition of the vector of translation~\cite{Ku}.  To make this precise, let us introduce the following notation: \bigskip

\noindent $\bullet \ $Let $\TT^d= \RR^d / \ZZ^d$ be the $d$-dimensional torus.

\noindent $\bullet \ $Let $\mu$ be the probability Haar measure on $\TT^d$.

\noindent $\bullet \ $The mapping $$T_{\alpha}: \TT^d \rightarrow \TT^d \qquad T_{\alpha}(x) = x + \alpha$$ is toral translation by vector $\alpha \in \RR^d$.

\noindent $\bullet \ $Denote the nearest distance of $\a \in \RR$ to the integers by $$\|\a\| := \inf_{p \in \ZZ}|\a -p|.$$  In higher dimensions, denote the nearest distance of $\a \in \RR^d$ to the integer lattice $\ZZ^d$ by $$\|\a\|_{\ZZ}:= \max_{i} \|\a_i\|,$$ where $\a_i$ is the $i$-th component of $\a$. 

\noindent $\bullet \ $For $d \geq 1$ and $\sigma \geq 0$, define
$$ \Omega_d(\sigma) := \{\alpha \in \RR^d \mid \exists C >
0 \textrm{ such that } \forall k \in \NN,  {\|k\alpha\|}_{\ZZ} \geq C
k^{-(1 + \sigma)/d}\}.$$

Let us denote $\Omega := \Omega_1$.  The sets, $\Omega_d(\sigma)$, are fundamental objects of study in the theory of Diophantine approximation.  In particular, $\Omega_d(0)$ is called the set of badly approximable vectors or vectors of constant type.  Note that $$\Omega_d(\sigma) \subseteq \Omega_d(\tau) \quad \textrm{ for } 0 \leq \sigma \leq \tau.$$  A more comprehensive introduction may be found in a number of sources, in particular~\cite{Sch}.

Kurzweil showed (and B. Fayad rediscovered~\cite{Fa}) the following for MSTP:

\begin{theo}[Kurzweil~\cite{Ku}]\label{theoKurzweil} The dynamical system $(\TT^d, \mu, T_{\alpha})$ has MSTP $\iff \alpha \in \Omega_d(0)$.
\end{theo}

\subsection{The Main Result}
Our main result is a generalization of Theorem~\ref{theoKurzweil} for $d=1$.  To begin, let us require sequences of measurable sets, $\{A_n\}_{n \in \NN}$, to satisfy a more stringent condition than (\ref{div}):  \begin{eqnarray} \label{div2} \sum_{n=1}^{\infty} (\mu (A_n))^s =  \infty \quad \textrm{ for }s \geq 1.\end{eqnarray}  We study an idea implicitly introduced by Fayad in his proof that mixing can occur in the absence of MSTP(~\cite{Fa},~Theorem~3.4):

\begin{defn}[$s$MSTP]
Let $s \geq 1$. The dynamical system $(M, \mu, T)$ has the \textbf{$s$-exponent monotone shrinking target property ($s$MSTP)} if, for any $x \in M$ and any $r \in {\mathcal DR}$ such that $A_n := B(x, r_n)$ satisfies (\ref{div2}), $\{A_n \}$ is BC for $T$.
\end{defn}  

\noindent Note that MSTP $\Rightarrow$ $s$MSTP.  \bigskip

The core of Fayad's argument in his Theorem~3.4\footnote{In~\cite{Fa}, there is a misprint in the statement of this theorem:  $\a \notin \Omega_d(2/d)$ should be replaced with $\a \notin \Omega_d(1)$.} is:  $\a \notin \Omega_d(1) \Rightarrow T_\a$ does not have $\frac{d+1} {d}$-MSTP.  He starts by constructing a certain flow on $\TT^{d+1}$ and noting that $\TT^d$ is a section on which the first return map is $T_\a$. Then he shows this core argument.  The sequence of monotonically shrinking balls on $\TT^d$ (constructed below in the proof of Proposition~\ref{propNotinNosMSTP}) which violates $\frac{d+1}d$-MSTP corresponds to a sequence of monotonically shrinking balls on $\TT^{d+1}$ whose sum of measures diverges, but whose limsup set has zero measure.  Thus, Fayad shows the time-$1$ map of the flow does not have MSTP.  For $d=3$, conjugating this map produces a real analytic map, which preserves the Haar measure and is mixing, but does not have MSTP(~\cite{Fa}, Corollary 3.5).

We show in Proposition~\ref{propNotinNosMSTP} that this core argument readily generalizes for all Diophantine exponents.  For $d=1$, the converse, $\alpha \in \Omega(s-1) \Rightarrow T_{\alpha}$ has $s$MSTP, is also true, but harder to prove than the comparable implication in~\cite{Ku} or~\cite{Fa}.  Our main result is this proof:

\begin{theo}\label{thmMR} Let $s \geq 1$.  The dynamical system $(\TT^1, \mu, T_\a)$ has $s$MSTP $ \iff \alpha \in \Omega(s-1).$
\end{theo}

\begin{proof} Apply Propositions~\ref{propNotinNosMSTP} and~\ref{propInsMSTP1} below. \end{proof}

Finally, it is worth noting that Theorem~\ref{thmMR} and~\cite{Fa} show that the tower of implications for the shrinking target properties $$\textrm{STP} \Rightarrow s\textrm{MSTP} \Rightarrow t\textrm{MSTP}$$ for $1 \leq s < t$ is strict.

\subsection{An Application.}  As an application of our technique from Section~\ref{secCRTHs}, we give in Section~\ref{secCRaA} a new, very short, proof of

\begin{theo} [Kim~\cite{Ki}]\label{thmKim} Let $\a \in \RR \backslash \QQ$.  Then $$\liminf n \|n \a - s\|_{\ZZ} = 0$$ for Lebesgue-a.e. $s \in \RR$.
\end{theo}

An introduction to Theorem~\ref{thmKim} can be found in~\cite{Ki}.  In this note, we consider this theorem from the alternate viewpoint of the shrinking target properties.  Consider the following admissible sets of radius sequences:  \begin{eqnarray} \label{A}\A_d := \{r \in \DR \mid \liminf n^{1/d}r_n > 0\}.\end{eqnarray}  By unraveling definitions, it is clear that, for all $\a \in \RR^d$, \begin{eqnarray} \label{KimSTP} \liminf n^{1/d} \|n \a -s\|_{\ZZ}=0 \textrm{ for a.e. } s \in \TT^d \iff T_\a \textrm{ has } \A_d\textrm{-STP}.\end{eqnarray}  Since a.e. $s \in \TT^d$ is equivalent to Lebesgue-a.e. $s \in \RR^d$, Theorem~\ref{thmKim} can be restated as: $$\textrm{Let } \a \in \RR \backslash \QQ. \textrm{  Then } (\TT^1, \mu, T_\a) \textrm{ has } \A_1\textrm{-STP}.$$

We immediately obtain: 

\begin{coro}[Logarithm Law for Circle Rotations] \label{corLogLaw} Let $\a \in \RR \backslash \QQ.$  Then $$\limsup \frac{-\log \|T^n_\a(x)\|_{\ZZ}} {\log n} =1$$ for a.e. $x \in \TT^1$.
\end{coro}

The above logarithm law is similar to Sullivan's logarithm law for geodesics found in~\cite{Su} and generalized in~\cite{KM}.  

\section{Toral Translations That Do Not Have $s$MSTP}\label{secTTTDNHs}

The following is a slight modification of the proof of Theorem 3.4 of~\cite{Fa}:

\begin{prop} \label{propNotinNosMSTP}
If $\alpha \notin \Omega_d(\sigma-1)$, then $T_{\alpha}$ does not have $\frac{d+\sigma-1}{d}$-MSTP.
\end{prop}

\begin{proof} Let us denote $s:=\frac{d+\sigma-1}{d}$.  We know that there exists $\{Q_n\} \subset \NN$ such that $Q_{n+1} \geq 2 Q_n$ and $\|Q_n \alpha\|_{\ZZ} \leq \frac {1}{2n^{2s+2/d} Q_n^{\sigma/d}}$.

Define $U_n := \lfloor n^{2s}Q_n^s \rfloor$ and $R_n := n^{-2/d}Q_n^{-1/d}$.  For $l \in [U_{n-1}, U_n-1]$, let $r_l :=R_n$.

Then we can easily show:  $\sum_{l=U_{n-1}}^{U_n-1} r_l^{ds} \geq (const) n^{2s} Q_n^s n^{-2s} Q_n^{-s} = const$. (Note that $const$ is a strictly positive constant.)  Summing over all intervals gives the divergence condition.

Also,

\begin{align*} \bigcup_{l=U_{n-1}}^{U_n-1} T_{\alpha}^{-l}B(x,r_l) &\subset \bigcup_{q=0}^{Q_n-1} \quad \bigcup_{p=0}^{\lceil n^{2s} Q_n^{s-1}\rceil} T_{\alpha}^{-pQ_n-q}B(x,R_n)\\  &\subset \bigcup_{q=0}^{Q_n-1} T_{\alpha}^{-q} B(x, 2R_n). \end{align*}

Thus, 

$$ \mu(\bigcup_{l=U_{n-1}}^{U_n-1} T_{\alpha}^{-l}B(x,r_l)) \leq Q_n \mu(B(x, 2R_n) = const/n^2.$$  Hence the limsup set has zero measure. \end{proof}

\section{Circle Rotations That Have $s$MSTP} \label{secCRTHs}

In this section, we prove:

\begin{prop} \label{propInsMSTP1} If $\alpha \in \Omega(s-1)$, then $T_{\alpha}$ has $s$MSTP.
\end{prop}

However, we will first give a simple proof of a weaker result:

\begin{prop} \label{propInsMSTP2} Let $ s > 1$ and $s-1 > \delta > 0 $ be arbitrary.  If $\alpha \in \Omega(s-1 - \delta)$, then $T_{\alpha}$ has $s$MSTP.
\end{prop}

In this section, $\alpha \in \RR \backslash \QQ.$  Without loss of generality, we may further assume that $ \alpha \in (-1/2, 1/2)$ as $T_{\alpha +1} = T_{\alpha}$.  Hence, $\alpha \in \RR \backslash \QQ \cap (-1/2, 1/2).$  Also, define $$x_n := T_{\alpha}^{n}(0) \quad \textrm{ for } n \in \NN \cup \{0\}.$$

We begin the proofs of Propositions~\ref{propInsMSTP1} and~\ref{propInsMSTP2} by noting some simple facts about continued fractions.  The terminology and standard results of the theory of continued fractions are found in many sources, in particular~\cite{Kh}.

\subsection{Continued Fractions} Let $p_i/q_i$ be the $i$-th order convergent of $\alpha$.
Define $$\Delta_i := |q_i \alpha-p_i| = \|x_{q_i}\|_{\ZZ}.$$  Since $\alpha$ is irrational, $\Delta_i > 0$ for all $i$.  We note that $q_0 =1$, $\Delta_0=|\alpha|$, and $\Delta_i < \frac 1 {2^{i/2}}$.  The following four lemmas are well-known:

\begin{lemma} \label{lemmMgap} For $i \geq 2$, $q_i \geq 2 q_{i-2}$.
\end{lemma}

\begin{lemma} \label{lemmUnif2} For all $i \in \NN$, $\frac 1 2 \Delta_{i-1}^{-1} < q_i \leq \Delta_{i-1}^{-1}$. 
\end{lemma}

\begin{lemm}
\label{lemmUnif}  
Let $0 \leq j < k <q_i$.  Then, $\|x_k-x_j\|_{\ZZ} \geq \Delta_{i-1}$.
\end{lemm}

\begin{lemm}[Diophantine Condition] \label{lemmDio} Let $\alpha \in \Omega(t -1).$  For all $i \in \NN$, $\Delta_i \geq C \Delta_{i-1}^{t}$.
\end{lemm}

We are now ready to apply these facts.

\subsection{Proof of the Propositions}  

Since all radii sequences under consideration are weakly monotonically decreasing, $\limsup B(x_n, r_n)$ is $T^{-1}_\a$-invariant.  To prove the Propositions, it suffices, by the ergodicity of $T^{-1}_\a$, to show that for a fixed $\eta > 0$ and all $r \in {\mathcal DR}$, there exists a $N_0$ such that $\mu(\cup_{n=0}^{N_0} B(x_n, r_n)) \geq \eta$.  

To begin the proofs, let us highlight a very simple, but key observation:

\begin{lemm} \label{propEZCase} Let $r \in {\mathcal DR}$ and let $0 < \varepsilon < 1/2$.  If there exists an $i \in \NN$, such that $r_{q_i-1} \geq \varepsilon \Delta_{i-1}$, then $\mu(\cup_{n=0}^{q_i-1} B(x_n, r_n)) \geq \varepsilon$.
\end{lemm}

\begin{proof} By Lemma~\ref{lemmUnif}, $\cup_{n=0}^{q_i-1} [x_n - \varepsilon \Delta_{i-1}, x_n + \varepsilon \Delta_{i-1}]$ is a disjoint union, which is contained in $\cup_{n=0}^{q_i-1} B(x_n,r_n)$.  Thus, by Lemma~\ref{lemmUnif2}, $$\mu(\cup_{n=0}^{q_i-1} B(x_n,r_n)) \geq \frac 1 2 \Delta_{i-1}^{-1} 2 \varepsilon \Delta_{i-1} = \varepsilon.$$ \end{proof}

\noindent{\bf Proof of Proposition~\ref{propInsMSTP2}}

\begin{proof} Let $ 0< C < 1$ be chosen so that $\|k \alpha\|_{\ZZ} \geq Ck^{-(s - \delta)}$ for all $k \in \NN$.  Let $0< \varepsilon < {(\frac C 4)}^{1/s}$ be chosen.  Let $\eta := \varepsilon/64$ ($\eta$ depends only on $\alpha$).

If, for all $i \in \NN$, $r_{q_i-1} < \frac \varepsilon {4} \Delta_{i-1}$, then by Lemma~\ref{lemmDio}, $r_{q_i-1}^s < \frac {\varepsilon^s} {4^s} \Delta_{i-1}^s \leq \frac {\varepsilon^s} {4^s C} \Delta_i \Delta_{i-1}^{\delta} \leq {\Delta_i \Delta_{i-1}^\delta}$.  Thus, $$ \sum_{i=1}^\infty \sum_{n = q_i}^{q_{i+1}-1} r_n^s \leq \sum_{i=1}^\infty q_{i+1} r_{q_i-1}^s \leq \sum_{i=1}^\infty \Delta_{i-1}^\delta \leq  \sum_{i=1}^\infty ((1/2)^{\delta/2})^i < \infty.$$  This is a contradiction, indicating that this possibility cannot happen.  

Otherwise, there exists an $i \in \NN$, such that $r_{q_i-1} \geq \frac \varepsilon {4} \Delta_{i-1}$.  Apply Lemma~\ref{propEZCase} to obtain $\mu(\cup_{n=0}^{q_i-1} B(x_n, r_n)) \geq \eta$. \end{proof} 

This pithy proof will not suffice for Proposition~\ref{propInsMSTP1}; modifications are necessary:\bigskip

\noindent{\bf Proof of Proposition~\ref{propInsMSTP1}}  
\begin{proof} Let $ 0< C < 1$ be chosen so that $\|k \alpha\|_{\ZZ} \geq Ck^{-s}$ for all $k \in \NN$.  Let $0< \varepsilon < {(\frac C 4)}^{1/s}$ be chosen.  Let $\eta := \varepsilon/64.$

\bigskip

\noindent {\bf Case 1.}  There exists an $i \in \NN$, such that $r_{q_i-1} \geq \frac \varepsilon {4} \Delta_{i-1}$.  \bigskip

Apply Lemma~\ref{propEZCase} to obtain $\mu(\cup_{n=0}^{q_i-1} B(x_n, r_n)) \geq \eta$. \bigskip

\noindent {\bf Case 2.}  For all $i \in \NN$, $r_{q_i-1} < \frac \varepsilon {4} \Delta_{i-1}$.  \bigskip
 
By Lemma~\ref{lemmMgap}, we obtain that for all $i \geq 0$, there exists $1 \leq j \leq 2$ such that $\lfloor \frac{q_{i+j}}{q_i} \rfloor \geq 2$ and $\forall \ 0 \leq k < j$, $\lfloor \frac{q_{i+k}}{q_i} \rfloor = 1$. Thus, $q_{i+j-1} < 2 q_i \Rightarrow \frac 1 2 \Delta^{-1}_{i+j-2} < 2 \Delta^{-1}_{i-1} \Rightarrow \Delta_{i-1} < 4 \Delta_{i+j-2}.$ Hence, we obtain a subsequence $\{q_{i_m}\}$ such that for all $i_m \leq k < i_{m+1}$, $\lfloor \frac{q_k} {q_{i_m}} \rfloor =1$, and $\lfloor \frac{q_{i_{m+1}}} {q_{i_m}} \rfloor \geq 2$.  We may start with $i_0=0$; hence $i_1=1$ as $q_0 =1$.  

For the remainder of the proof of this proposition, $m$ is a strictly positive integer.  Define $K_m:= \lfloor \frac{q_{i_{m+1}}} {q_{i_m}} \rfloor$ and note that $\Delta_{i_m-1} < 4 \Delta_{i_{m+1}-2}.$

\begin{sublemm} \label{sublemmsmallrs}For $q_{i_m} \leq N < q_{i_{m+1}}$, one has $r_N^s < \frac{\Delta_{i_{m+1}-1}}{4}.$  Moreover, the balls $\{B(x_N, r_N^s)\}_{N=q_{i_m}}^{q_{i_{m+1}}-1}$ are pairwise disjoint.
\end{sublemm}

\begin{proof} By Lemma~\ref{lemmDio}, $r_{q_{i_m}-1} < \frac \varepsilon {4} \Delta_{i_m-1} < \varepsilon \Delta_{i_{m+1}-2} \Rightarrow r_{q_{i_m}-1}^s < \varepsilon^s  \Delta^s_{i_{m+1}-2} \leq \frac {\varepsilon^s} {C} \Delta_{i_{m+1}-1} < \frac{\Delta_{i_{m+1}-1}}{4}.$  Pairwise disjointness follows from Lemma~\ref{lemmUnif}. \end{proof}

\begin{sublemm} \label{sublemmdisjoint1} Let $0 \leq n < q_{i_m} \leq N <q_{i_{m+1}}$.  If $B(x_n, r^s_n) \cap B(x_N, r^s_N) \neq \emptyset$, then $r_n \geq \frac {\varepsilon} {4} \Delta_{i_m-1}$.
\end{sublemm}

\begin{proof} By Lemma~\ref{lemmUnif} and Sublemma~\ref{sublemmsmallrs}, $r^s_n \geq \Delta_{i_{m+1}-1} - \frac{\Delta_{i_{m+1}-1}}{4}$. Thus by Lemma~\ref{lemmDio}, $r_n \geq (3/4)^{1/s} C^{1/s} \Delta_{i_{m+1}-2} \geq 3^{1/s} \varepsilon \Delta_{i_{m+1}-2} \geq \frac \varepsilon  4 \Delta_{i_{m}-1}$. \end{proof} 

Let $1 \leq j \leq K_m-1$.  Define: $$A_m := \bigcup_{n=0}^{q_{i_m}-1}B(x_n, r_n)$$ $$\D_m := \{B(x_n, r_n^s)\}_{n=q_{i_{m-1}}}^{q_{i_m}-1} \quad D_m := \coprod_{n=q_{i_{m-1}}}^{q_{i_m}-1} B(x_n, r_n^s)$$ $$\hspace{10.5mm} \B_{jm} := \{B(x_n, r_n^s)\}_{n=j q_{i_m}}^{(j+1)q_{i_m}-1}.$$ 

\bigskip

We wish to refine Sublemma~\ref{sublemmdisjoint1}:

\begin{sublemm} \label{sublemmdisjoint2}  Let $0 \leq n < q_{i_m}$.  If $B(x_n, r^s_n)$ intersects $L>1$ elements of $\B_{jm}$, then $r^s_n \geq \frac{L-1}{4} \Delta_{i_m-1}$.
\end{sublemm} 

\begin{proof} By Lemma~\ref{lemmUnif}, the least distance between centers of any two distinct elements of $\B_{jm}$ is at least $\Delta_{i_m-1}$.  Thus, by Sublemma~\ref{sublemmsmallrs}, we obtain:  $r^s_n \geq  \frac{L-1} {2} \Delta_{i_m-1} - \frac{\Delta_{i_{m+1}-1}}{4}.$ \end{proof}

\bigskip

\noindent{\bf A Special Case} \bigskip

Before turning to the full generality of Case 2, let us prove another special case:  

\begin{sublemm}\label{sublemmEqualRadPerCycle}  Let $r_n = R_m$ for $q_{i_{m-1}} \leq n < q_{i_m}$.  Then there exists an $N_0 \in \NN$ such that $\mu(\cup_{n=0}^{N_0} B(x_n, r_n)) \geq \eta$.
\end{sublemm}

\noindent{\bf Remark:}  Of course, $\sum_n r_n^s = \infty$ and, for all $i \in \NN$, $r_{q_i-1} < \frac \varepsilon {4} \Delta_{i-1}$ still hold.  \bigskip  

The following proof of Sublemma~\ref{sublemmEqualRadPerCycle} is similar to a proof in~\cite{Fa}.

\begin{proof} To begin, we observe:

\begin{sublemm} \label{sublemmdisjoint3} The balls $\{B(x_N, \frac {\Delta_{i_{m+1}-1}}{2})\}_{N=q_{i_m}}^{q_{i_{m+1}}-1}$ are pairwise disjoint.
\end{sublemm}

\begin{proof} Apply Lemma~\ref{lemmUnif}.  \end{proof}  

Now we need another refinement of Sublemma~\ref{sublemmdisjoint1}:

\begin{sublemm} \label{sublemmdisjoint4} Let $q_{i_m} \leq N <q_{i_{m+1}}$.  If $\cup_{p=1}^m D_p \cap B(x_N, r^s_N) \neq \emptyset$, then $\cup_{p=1}^m D_p \supset (x_N - \frac {\Delta_{i_{m+1}-1}}{2}, x_N - \frac {\Delta_{i_{m+1}-1}}{4})$ or $\cup_{p=1}^m D_p \supset (x_N + \frac {\Delta_{i_{m+1}-1}}{4}, x_N + \frac {\Delta_{i_{m+1}-1}}{2})$.
\end{sublemm}

\begin{proof} Follows from the first line of the proof of Sublemma~\ref{sublemmdisjoint1}.  \end{proof}  

A corollary of Sublemmas~\ref{sublemmdisjoint3} and~\ref{sublemmdisjoint4} is:

\begin{sublemm}\label{sublemmdisjoint5} Let $P$ elements of $\D_{m+1}$ intersect $\cup_{p=1}^m D_p$.  Then, $\mu(\cup_{p=1}^m D_p) \geq \frac {P} {4} \Delta_{i_{m+1}-1}$.
\end{sublemm}

\begin{proof} Let $x_{N_1}, \cdots, x_{N_P}$ be the centers of these $P$ elements of $\D_{m+1}$.  By Sublemma~\ref{sublemmdisjoint4}, $\cup_{p=1}^m D_p$ contains a piece of the open ball $B(x_{N_k}, \frac {\Delta_{i_{m+1}-1}}{2})$ of measure at least $\frac { \Delta_{i_{m+1}-1}} {4}$ for every $1 \leq k \leq P$.  By Sublemma~\ref{sublemmdisjoint3}, there are no overlaps.  Hence the result follows.  \end{proof}  

If $\lfloor \frac {q_{i_{m+1}} - q_{i_m}} 2 \rfloor$ elements of $\D_{m+1}$ intersect $\cup_{p=1}^m D_p$, then, $$\mu(\cup_{p=1}^m D_p) \geq \frac { \lfloor \frac {q_{i_{m+1}} - q_{i_m}} 2 \rfloor} {4} \Delta_{i_{m+1}-1}.$$  Now $K_m \geq 2 \Rightarrow q_{i_m} \leq q_{i_{m+1}}/2$. Hence, $$\mu(\cup_{p=1}^m D_p) \geq \frac { \lfloor \frac {q_{i_{m+1}}} {4}\rfloor} {4} \Delta_{i_{m+1}-1} \geq \frac 1 {64} \Delta_{i_{m+1}-1}^{-1} \Delta_{i_{m+1}-1} > \eta,$$ and we have shown Sublemma~\ref{sublemmEqualRadPerCycle} for this possibility.

Otherwise, at least $\lceil \frac {q_{i_{m+1}} - q_{i_m}} 2 \rceil$ elements of $\D_{m+1}$ are disjoint from $\cup_{p=1}^m D_p$.  By Sublemma~\ref{sublemmsmallrs} , these elements are also pairwise disjoint; hence we obtain:  \begin{eqnarray} \label{eqnEqualRadPerCycle} \mu(\cup_{p=1}^{m+1} D_p) \geq \mu(\cup_{p=1}^m D_p) + \frac {q_{i_{m+1}} - q_{i_m}} 2 \mu(B(x_{q_{i_{m+1}}-1}, R_{m+1}^s)). \quad \end{eqnarray}

This process is recursive.  If the recursion stops, then we have shown the Sublemma; otherwise, we obtain by (\ref{eqnEqualRadPerCycle}):  $$ \mu(\cup_{m=1}^\infty D_m) \geq \sum_{m=1}^\infty (q_{i_{m+1}} - q_{i_m}) R_{m+1}^s = \sum_{n=q_1}^\infty r_n^s = \infty.$$

This is a contradiction, indicating the recursion must stop and, thus, proving Sublemma~\ref{sublemmEqualRadPerCycle}. \end{proof} \bigskip

\noindent{\bf The General Proof} \bigskip

Let us consider the general case.  We wish to recursively define a subset of $\cup_{m=1}^\infty D_m$, call it $E_\infty$, as follows.  Let $\E_2 := \D_2$ and $E_2 := D_2$.  (Note that $E_2$ is a disjoint union by Sublemma~\ref{sublemmsmallrs}.)

\bigskip

\noindent{\bf The recursive algorithm follows:}

\bigskip

\noindent{\bf Enough balls intersect}  \bigskip

Assume that there are at least $\lfloor q_{i_2}/2 \rfloor$ elements of $\B_{j2}$ that intersect $E_2$; call this set $\hat{\B}_{j2}$.  (Note that $q_{i_2} \geq 4$.)  Let us count the elements of $\hat{\B}_{j2}$ as follows.  There is an element (not necessarily unique) of $\E_2$ that intersects the most number of elements of $\hat{\B}_{j2}$.  Call this number $L_1$.  Remove these $L_1$ elements from $\hat{\B}_{j2}$.  Repeat.  Since $card(\hat{\B}_{j2}) < \infty$, this process terminates.  Let $L_1 \geq L_2 \geq \cdots \geq L_N$ be the resulting sequence. Let $B(x_{n_1}, r_{n_1}^s), \cdots, B(x_{n_N}, r_{n_N}^s)$ be the corresponding elements of $\E_2$.  Note that $$L_1 + \cdots + L_N = card(\hat{\B}_{j2}) \geq \lfloor q_{i_2}/2 \rfloor.$$ Also, note that $B(x_{n_h}, r_{n_h}^s)$ intersects at least $L_h$ elements of $\hat{\B}_{j2}$.  It may be possible for it to intersect more.

Let $N \geq H \geq 0$ be the greatest such that $L_H > 1$.  Then $$L_{H+1} + \cdots + L_{N} = N -H$$ since it is a sum of $1$'s (or possibly zero).  Also, $H \leq card(\hat{\B}_{j2})/2$.  There are two cases:

Let $N-H \leq \lfloor q_{i_2}/2 \rfloor/4$.  (Note that if $H=0$, $N = card(\hat{\B}_{j2})$.  Hence, $N-H > \lfloor q_{i_2}/2 \rfloor/4$.)  By Sublemma~\ref{sublemmdisjoint2}, for $1 \leq h \leq H$, $r_{n_h}^s \geq \frac {L_h-1}{4} \Delta_{i_2-1}$.  The elements in $\E_2$ are disjoint; thus, \begin{align*}\mu(E_2) &\geq \sum_{h=1}^{H} \mu(B(x_{n_h}, r_{n_h}^s) \geq \frac 1 2 (L_1 + \cdots + L_H -H) \Delta_{i_2-1}\\ &\geq \frac 1 2 (card(\hat{\B}_{j2})/2-(N-H))\Delta_{i_2-1} \geq \frac 1 {8} \lfloor q_{i_2}/2 \rfloor \Delta_{i_2-1}\geq \frac 1 {64} > \eta.\end{align*}

Otherwise, $N-H > \lfloor q_{i_2}/2 \rfloor/4$.  Then the elements of $\E_2$, \newline $B(x_{n_{H+1}}, r_{n_{H+1}}^s), \cdots, B(x_{n_N}, r_{n_N}^s)$, intersect at least one element of $\hat{\B}_{j2}$.  Consequently, Sublemma~\ref{sublemmdisjoint1} implies that, for $H+1 \leq h \leq N$, one has $r_{n_h} \geq \frac {\varepsilon} {4} \Delta_{i_2-1}$.  Thus, by Lemma~\ref{lemmUnif}, $$A_2 \supset \cup_{h=H+1}^N B(x_{n_h}, r_{n_h}) \supset \coprod_{h=H+1}^N B(x_{n_h}, \frac {\varepsilon} {4} \Delta_{i_2-1}).$$  Hence, $\mu(A_2) \geq \frac{\varepsilon}{8} \lfloor q_{i_2}/2 \rfloor \Delta_{i_2-1} \geq \frac \varepsilon {64} \geq \eta$.

\bigskip

\noindent{\bf Not enough balls intersect} \bigskip

Otherwise, at least $\lceil q_{i_2}/2 \rceil$ elements of $\B_{j2}$ are disjoint from $E_2$ for each $1 \leq j \leq K_2-1$.  Call this subset of $\B_{j2}$, $\tilde{\B}_{j2}$. Then define $\E_3 := \E_2 \cup_{j=1}^{K_2-1} \tilde{\B}_{j2}$, and $E_3 := \cup_{B \in \E_3} B$. 
By Sublemma~\ref{sublemmsmallrs}, $E_3$ is also a disjoint union of balls. \bigskip

\noindent{\bf Repeat the recursion.} \bigskip

If the recursion stops, then we have shown the Proposition.  Otherwise, we obtain $E_\infty = D_2 \coprod \coprod_{m=2}^\infty \coprod_{j=1}^{K_m-1} \tilde{\B}_{jm}$ as a subset of $\cup_{m=1}^\infty D_m$.  Hence, \begin{eqnarray} \label{eqnunionsum} \mu(E_\infty) \geq 2 \sum_{m=1}^\infty \sum_{j=1}^{K_m-1} \lceil q_{i_m}/2 \rceil r_{(j+1)q_{i_m}}^s \ .\end{eqnarray}

Now $$\infty = \sum_n r_n^s \leq \sum_{m=1}^{\infty} [\sum_{j=1}^{K_m-1} q_{i_m} r_{jq_{i_m}}^s + \sum_{l = K_m q_{i_m}}^{q_{i_{m+1}-1}} r_l^s] \leq \sum_{m=1}^{\infty} \sum_{j=1}^{K_m-1} 2 q_{i_m} r_{jq_{i_m}}^s.$$  The last inequality follows from the fact that $\sum_{l = K_m q_{i_m}}^{q_{i_{m+1}-1}} r_l^s$ has fewer than $q_{i_m}$ terms. Hence, $\sum_{m=1}^{\infty} \sum_{j=1}^{K_m-1} \lceil q_{i_m}/2 \rceil r_{jq_{i_m}}^s = \infty.$

Now the internal sum is telescoping:  $$0 \leq \sum_{j=1}^{K_m-1} \lceil q_{i_m}/2 \rceil r_{jq_{i_m}}^s - \sum_{j=1}^{K_m-1} \lceil q_{i_m}/2 \rceil r_{(j+1)q_{i_m}}^s \leq \lceil q_{i_m}/2 \rceil r_{q_{i_m}}^s .$$

Hence, $$0 \leq \sum_{m=1}^{\infty} \sum_{j=1}^{K_m-1} \lceil q_{i_m}/2 \rceil r_{jq_{i_m}}^s - \sum_{m=1}^{\infty} \sum_{j=1}^{K_m-1} \lceil q_{i_m}/2 \rceil r_{(j+1)q_{i_m}}^s \leq \sum_{m=1}^{\infty} \lceil q_{i_m}/2 \rceil r_{q_{i_m}}^s.$$

We now obtain two cases: \bigskip

\noindent{\bf Case 2A:}  $\sum_{m=1}^{\infty} q_{i_m} r_{q_{i_m}}^s < \infty.$  \bigskip

This case implies that the right hand side of (\ref{eqnunionsum}) diverges, a contradiction indicating that, for Case 2A, the recursion must stop, and, hence, proving the Proposition.  \bigskip

\noindent{\bf Case 2B:}  $\sum_{m=1}^{\infty} q_{i_m} r_{q_{i_m}}^s = \infty.$ \bigskip

By the definition of the subsequence $i_m$, $q_{i_m} \geq 2 q_{i_{m-1}}$.  Thus, $\sum_{m=1}^{\infty} (q_{i_m}-q_{i_{m-1}}) r_{q_{i_m}}^s = \infty$.  Let us now define $r' \in {\mathcal DR}$ by $r_n' := r_{q_{i_m}}$ for $n = q_{i_{m-1}}, \cdots, q_{i_m}-1$.  Hence, $r_n' \leq r_n$.  We wish to prove the Proposition for this new set of balls.  Since the new radii are smaller, the condition of Case 2 applies.  Applying Sublemma~\ref{sublemmEqualRadPerCycle} completes the proof of the Proposition.  \end{proof}

\section{Circle Rotations and $\A$STP} \label{secCRaA}

In this section, we give a new, very short, proof of Theorem~\ref{thmKim} by applying the very simple observation about continued fractions that is shown in Lemma~\ref{propEZCase}.  First, recall the definition of $\A_1$ from (\ref{A}).\bigskip

\noindent{\bf Proof of Theorem~\ref{thmKim}}  
\begin{proof} Let $\tilde{r} \in \A_1$.  Then for every $1 > \delta > 0$ small enough, there exists $N \in \NN$ such that, for all $n \geq N$, $\tilde{r}_n \geq \frac \delta n$.  Let $k$ be a strictly positive integer.  Define $r_n := \frac \delta {n +k}$.  Then there exists an $i \in \NN$ such that $q_i > \max(k,N)$.  Hence, $r_{q_i-1} \geq \frac \delta {2 (q_i-1)} \geq \frac \delta 2 \Delta_{i-1}$.  Applying Lemma~\ref{propEZCase} and $T_{-k \a}$, $$\mu(\cup_{n=k}^{\infty} B(x_n,\frac \delta n)) = \mu(\cup_{n=0}^{\infty} B(x_{n+k},r_n))= \mu(\cup_{n=0}^{\infty} B(x_n,r_n)) \geq \delta/2.$$

Hence, $\mu(\cap_{k=N}^\infty \cup_{n=k}^{\infty} B(x_n,\frac \delta n)) = 1$.  Thus, $T_\a$ has $\A_1$-STP, which, by (\ref{KimSTP}), is equivalent to the conclusion of Theorem~\ref{thmKim}. \end{proof}

\section{Open Questions}

The most immediate questions concern generalizations to higher dimensions, namely to translations on $\TT^d$ for $d >1$.  For generalizing Theorem~\ref{thmMR}, we know that, for $s \geq 1$, $\alpha \in \Omega_d(0) \Rightarrow T_{\alpha}$ has $s$MSTP by Theorem~\ref{theoKurzweil}, and, in the converse direction, that, for $s \geq 1$, $\alpha \notin \Omega_d(sd-d) \Rightarrow T_{\alpha}$ does not have $s$MSTP by Proposition~\ref{propNotinNosMSTP}.  What condition in higher dimensions yields an ``if and only if'' as it is in Theorem~\ref{thmMR} for $d = 1$?

For generalizing Theorem~\ref{thmKim} or the logarithm law, Theorem~\ref{theoKurzweil} and the (doubly metric) inhomogeneous Khinchin-Groshev theorem (see ~\cite{Ca}, Chapter VII, Theorem II) imply that,  for all $\a \in \Omega_d(0)$ and almost all $\a \in \RR^d$ respectively, \begin{eqnarray}\label{eqnkimHD}\liminf n^{1/d} \|n \a -s\|_{\ZZ}=0 \textrm{ for Lebesgue-a.e. } s \in \RR^d,\end{eqnarray} and \begin{eqnarray} \label{loglawHD} \limsup \frac{-\log \|T^n_\a(x)\|_{\ZZ}} {\log n} = \frac 1 d \quad \textrm{ for a.e. } x \in \TT^d.\end{eqnarray}  As noted before, (\ref{eqnkimHD}) is equivalent to $$(\TT^d, \mu, T_\a) \textrm{ has } \A_d\textrm{-STP.}$$  Is it true, as it is in Theorem~\ref{thmKim}, that all minimal $T_\a$ on $\TT^d$ have $\A_d$-STP?  Is it true, as it is in Corollary~\ref{corLogLaw}, that all minimal $T_\a$ on $\TT^d$ satisfy the logarithm law (\ref{loglawHD})? 

\section{Another Shrinking Target Property}

After this note was submitted, Theorem 3 of~\cite{Ku} was brought to the author's attention. In this section, we compare this theorem and Theorem~\ref{thmMR}.

Consider the following sets of admissible radius sequences \[\B_s := \{r \in \DR \mid \{r_n\} \textrm{ satisfies condition ($\ast$)}\},\] where  condition ($\ast$), found in~\cite{Ku}, is defined as follows:  \begin{itemize} \item[($\ast$)] There exists a function $\delta(q)$ defined for $q\geq1$, $\delta(q)\geq1$, $\delta(q) \rightarrow \infty$ steadily as $q \rightarrow \infty$ and a sequence of positive integers $t_1 < t_2 < t_3 < \cdots$, \[ t_{i+1} > t_i^s (\delta(t_i))^{s+1},\] such that  \[\sum_{i=1}^\infty t_i r_{\lfloor t_i^s (\delta(t_i))^{s+1} \rfloor} = \infty.\] \end{itemize}

For comparison with Theorem~\ref{thmMR}, fix $\varphi(q)= q^{-(s+1)} $ in Theorem 3 of~\cite{Ku} and translate this special case of the theorem into the viewpoint of the shrinking target properties:   \begin{itemize} \item[] Let $s \geq 1$.  The dynamical system $(\TT^1, \mu, T_\a)$ has $\B_s$-STP $ \iff \alpha \in \Omega(s-1).$ \end{itemize}

On the other hand, for $(\TT^1, \mu, T_\a)$, $s$MSTP is equivalent to $\DR_s$-STP where \[\DR_s := \{r \in \DR \mid \sum_n r^s_n = \infty\}.\]  It is known that $\B_1 = \DR_1$~\cite{Ku}.  Neither inclusion, however, is true for $s > 1$.

\begin{prop} \label{thmKsMSTP}Let $s > 1$.  Define $r_n = \frac 1 {(n \log n)^{1/s}}$.  Then $r \in \DR_s \backslash \B_s$.
\end{prop}

\begin{proof} Clearly, $r \in \DR_s$.  Assume that $r \in \B_s$.  Then there exists $\{t_n\}$ and $\delta$ as in condition ($\ast$).  Now $t_1 \geq 1$ and, for $n \geq 2$,  $t_n \geq 2^{s^{n-2}}.$  Thus, for $N$ large enough, \[ \infty = \sum_{n=N}^\infty t_n r_{\lfloor t_n^s (\delta(t_n))^{s+1} \rfloor} \leq \sum_{n=N}^\infty \frac{t_n} {(t_n^s \log(t_n^s))^{1/s}} \leq \frac 1 {\log(2)^{1/s}}\sum_{n=N-1}^\infty \frac 1 {s^{n/s}} < \infty, \] a contradiction.
\end{proof}

\begin{remark} Many other radius sequences such as $\{\frac 1 {(n \sum_{k=1}^n \frac 1 k)^{\frac 1 s}}\}$, $\{\frac 1 {(n \log n \log \log n)^{\frac 1 s}}\}$, etc., are in $\DR_s \backslash \B_s$ for $s > 1$.\end{remark}

\begin{prop} For $s > 1$, $\B_s \backslash \DR_s \neq \emptyset$.
\end{prop}

\begin{proof}  Choose $\{t_n\}$ as in condition ($\ast$) and define $\delta(t_n) = n^{\frac {s-1} {2(s+1)}}$.  For \[\lfloor t_{n-1}^s (\delta(t_{n-1}))^{s+1} \rfloor < n \leq \lfloor t_n^s (\delta(t_n))^{s+1} \rfloor,\] define $r_n = \frac 1 {t_n n}$.  Thus, $r \in \B_s$.  

But, \[\sum_{n=1}^\infty r^s_n \leq \sum_{n=1}^\infty t_n^s (\delta(t_n))^{s+1} r_{\lfloor t_n^s (\delta(t_n))^{s+1} \rfloor}^s  = \sum_{n=1}^\infty \frac 1 {n^{\frac{s+1} 2}} < \infty. \qedhere\] \end{proof}

Hence, neither Theorem~\ref{thmMR} nor Theorem 3 of~\cite{Ku} is a consequence of the other.  Thus, if we define  ${\mathcal R}_s = \DR_s \cup \B_s,$ the two theorems together imply a stronger theorem:
 
\begin{theo}\label{thmMRKur}Let $s \geq 1$.  The dynamical system $(\TT^1, \mu, T_\a)$ has ${\mathcal R}_s$-STP $ \iff \alpha \in \Omega(s-1).$ \end{theo}

\begin{proof} Apply Theorem 3 of~\cite{Ku} and Theorem~\ref{thmMR}. \end{proof}

Is there a set $\A$ strictly larger than ${\mathcal R}_s$ such that, if we replace ${\mathcal R}_s$ with $\A$ in Theorem~\ref{thmMRKur}, the theorem still holds?  How large can such an $\A$ be?

\section*{Acknowledgements} The author would like to thank his advisor, Professor Dmitry Kleinbock, for many helpful discussions and the referee for many helpful suggestions.  The author would also like to thank Dr. Alan Haynes for bringing Theorem 3 of~\cite{Ku} to the author's attention.

\end{document}